\documentclass[reprint, a4paper, nofootinbib, showpacs,
showkeys, twoside]{revtex4-1}
\usepackage[cp1251]{inputenc}
\usepackage[T2A]{fontenc}
\usepackage[greek,english]{babel}
\usepackage{amsthm,amsmath,graphicx,color}
\usepackage[psamsfonts]{amssymb}
\usepackage{mathrsfs}

\usepackage[colorlinks,bookmarks=false,pagebackref=false]{hyperref}
\hypersetup{pdfauthor={Yu. Brezhnev}, linktocpage=false,
		pdfwindowui=false, pdfmenubar=false, 
		citecolor=blue, pdfstartview=FitH, bookmarksopen=false}

\usepackage{YURA,Ulem}
\bodycenter[13mm]{0mm}{0mm} 

\renewcommand{\,}[1][1]{%
	\ifmmode\mkern#1mu\else\kern#1\dimexpr0.05556em\fi\relax}%
\renewcommand{\!}[1][1]{\,[-#1]}
\renewcommand{\.}{\,[3.]}

\def\yblue{\blue}

\def\+{  \mathbin{\Over[1.2ex]{\scalebox{0.5}[1]{$\sss\bo\frown$}}
	{\vcenter{\hbox{\small$+$}}}}}
\def\PM+{\mathbin{\Over[1.4ex]{\scalebox{0.5}[1]{$\sss\bo\frown$}}
	{\vcenter{\hbox{\small$\pm$}}}}}

\renewcommand{\|}[1][0]{\vcenter{\hbox{%
	\ifcase#1$|\hspace{-0.18em}|$%
	\or\small$\big|\hspace{-0.23em}\big|$%
	\else${\left|\vbox to5.3ex{}\hspace{-0.23em}\right|}$%
	\fi}}\relax}

\newcommand{\mbig}[2][4]{\vcenter{\hbox{%
	\ifcase#1$#2$\or\small$\big#2$\else%
	${\left#2\vbox to5.3ex{}\right.}$\fi}}\relax}

\makeatletter
\def\@keys@name{\emph{Key words}: }%
\makeatother

\def\tgreek#1{\textgreek{\texttt{\upshape #1}}} %
\def\greekA{\skew{-2}\vec{\tgreek{a}}}
\def\greekB{\skew{-3}\vec{\tgreek{b}}}
\def\greekC{\skew{-3.5}\vec{\tgreek{g}}}
\def\greekD{\skew{-3}\vec{\tgreek{d}}}
\def\greekAo{\skew{-1}\vec{\tgreek{a}}_%
	{\raisebox{0.2ex}{$\sss\,\msf o$}}}
\def\greekBo{\skew{-2}\vec{\tgreek{b}}_%
	{\raisebox{0.2ex}{$\sss\msf o$}}}

\def\vecO{\skew{-3}%
	\vec{\hbox{\footnotesize$\vbox to1.4ex{}\smash{\msf0}$}}}
\def\boB {{\ds\bo{\frak B}^{\![3]\sp}}}
\def\boA {{\ds\bo{\frak A}}}
\def\boAo{{\ds\bo{\frak A}^{\!\sss\msf o}}}
\def\VV{V{\Times\,}V}

\newcommand{\EqDown}[3][3.3ex]{\mathbin{\Under[#1]
	{\makebox[0ex][r]{#2}{\ds\Downarrow}\makebox[0ex][l]{#3}}{\ds=}}}
\newcommand{\LongEqDown}[3][3.3ex]{\Under[#1]{
	\makebox[0.5em][r]{#2\hskip0.59em}%
	\rotatebox[origin=cc]{90}{$\smash\Longleftarrow$}
	\makebox[0ex][l]{#3}}{\ds{}={}}}

\begin{document}
\title{A Note on the Measure of Vector and Pythagorean Theorem}
\author{\vskip-1exYurii V.~Brezhnev}
\affiliation{\vskip1ex Immanuel Kant Federal University,
K\"onigsberg, Russia}
%



\begin{abstract}\sffamily
Why the square? We present a geometry"=axiom"=free derivation of the
Pythagorean theorem and the square at its core, establishing their
algebraic origin from within the bare vector"=space framework. Such
concepts as the (right) angle, rotation, inner product,
orthogonality \etc\ also emerge as a logical construct rather than
taken as given. They are necessitated by the square, and the ensuing
theory, in turn, \emph{canonically} stems from a \emph{single}
definitional primitive""---the ($\ds\bbR^{\!\sp}\![3]$"=quantitative)
invariant $\cal Q$"=measure of a vector. This provides the core of an
algebraic justification for Euclidean geometry. Equally important,
these findings account (also canonically) for the complex
mod"=squared $p = |\frak a|^2$\,---the quantum Born rule---and point
out what is even admissible for being quantitatively interpreted.
The linear structure and its endomorphisms are rigid in the sense
that the \emph{well"=defined} interpretable turns out to be, up to
gauge $\cal Q\goto[1]\const{\,\times}\cal Q$, the unique
gauge"=invariant measure $\cal Q=\|{\cdot}{\cdot}{\cdot}\|^2$;
independently of the field $\bbR$ or $\bbC$.
\end{abstract}

\maketitle

{
\def\hboxDots{\leaders\hbox{\footnotesize.\kern\fontdimen3\font}}
\noindent\hyperlink{I1}{\bfcour{\yblue\small Introduction}}%
\hboxDots\hfill1

\noindent\hyperlink{I2}
{\bfcour{\yblue\small Quantification of vector}}%
\hboxDots\hfill1

\noindent\hyperlink{Inew}
{\bfcour{\yblue\small Extending the $V\![4]$-space algebra:\
$\cal Q$-measure}}%
\hboxDots\hfill2

\noindent\hyperlink{I3}
{\bfcour{\yblue\small Definition $\goto$ Theorem\,\,}}%
\hboxDots\hfill3

\noindent\hyperlink{I4}
{\bfcour{\yblue\small How~\,\,does~\,\,Pythagorean~\,\,square~\,\,%
come~\,\,into~\,\,being?}}%
\hboxDots\hfill4

\noindent\hyperlink{I5}
{\bfcour{\yblue\small Theorem $\goto$ Definitions}}%
\hboxDots\hfill5

\noindent\hyperlink{I6}
{\bfcour{\yblue\small The main theorem:\ canonicity}}%
\hboxDots\hfill5

\noindent\hyperlink{I7}
{\bfcour{\yblue\small References}}\hboxDots\hfill6%
}%
\hypertarget{I1}{}

\medskip
\noindent%
\bfcour{\yblue Introduction.}
The vector"=space ($V\!\!$"=space) formulation of the Pythagorean
theorem

\vbox to4.2ex{
\begin{equation}\label{pif}
\|\greekA \+ \greekB\|^2 = \|\greekA\|^2 + \|\greekB\|^2
\qquad(\greekA \mathrel\bot \greekB)
\end{equation}}
\vspace{1.1ex}
\noindent
and the square that lies at its core are a common topic in textbooks
\cite{halm}. However the square per~se is either postulated, as in
normed spaces, or is a consequence of definition via a binary
structure on the space---the scalar"=product $\langle\greekA\,\,{,}
\greekA\rangle=\|\greekA\|^2\!\!$, supplemented with the
parallelogram~law~$\|\greekA {\,\,\+\,\,} \greekB\|^2\!+\,
\|\greekA\mathbin{\,\Over[1ex]
{\scalebox{0.5}[1]{$\sss\bo\frown$}}{\raisebox{0.22pt}[1.12ex]
	{\hbox{\small$-$}}}\,\,}\greekB\|^2 {\,\,=\.} 2\,\|\greekA\|^2
{\,+\,\,} 2\,\|\greekB\|^2$ for a norm~$\|{\cdot}{\cdot}{\cdot}\|$.

Although the theorem is not just old but ancient \cite{maor}, the
questions surrounding such a norm still persist in contemporary
literature: ``So, why squares?''\@ \cite[p.~262]{givental}. Why not
cube or some other $\ds\ell^\textsf{p}_{\mathstrut}\!$"=norm \cite{aaronson}?
Whence such quadratic expressions? What, then, is the status of
orthogonality? How are we to account for it? Do polarization
identities \cite{halm}, the Riesz--Fr\'echet representation theorem,
and other ingredients of the theory lie at its foundation and are
necessary for it? Moreover, the same considerations are of great
importance for quantum theory. See Refs.~\cite{aaronson, br2} for
discussion.

Presented as a brief communication, this note aims to demonstrate
the algebraic nature of the Pythagorean theorem and the resulting
geometry. Thereby, they become not foundational givens but emergent
properties of an algebraic framework. Stated differently, we
establish that Pythagorean two, the concept of a right angle
\cite{maor}, and other related constructs do not require ad~hoc
postulation and are the \emph{structural necessity} of the
$V\!\!$"=space algebra and its endomorphisms; \ie, they are the
\emph{canonical extensions}. More
precisely, one can state the following.%
\hypertarget{P}{}

\vskip0.4ex
\textbf{Proposition.} \emph{Pythagoras\!'\,\ theorem, along with its
geometrical premises and consequences, stems from a single,
axiom-free and natural definition. The definition introduces a
numerical quantity that is invariantly associated with vector;
quantifying the vector by a $\cal Q$"=measure.}

To the best of our knowledge, such a low"=level re"=examination of the
grounds for the theorem and its status (theorem \goto\
\emph{self"=suggested} definition, followed by hierarchy of concepts)
appears to have escaped detailed scrutiny in the literature
\cite{maor}. The meaning of the terms `invariantly, quantifying by a
$\cal Q$"=measure, single definition', as well as the requisite
technicalities, will be described in the subsequent sections.

Taken together, these findings also show the strong similarity,
almost identity, between the Pythagorean square over $\bbR$ and the
quantum Born rule over $\bbC$ \cite{aaronson, br2}. The mechanism
behind the emergence of the complex modulus"=squared $p =
|\frak a|^2$ is a long"=standing challenge in quantum
foundations~\cite{aaronson}. As we shall see in
\hyperlink{D3}{\black Corollaries~\brown\textbf{2}} and
\hyperlink{C3}{\black \brown\textbf{3}}, the resolution of this issue
is mathematical in character rather than physical.

\hypertarget{I2}{}

\vskip0.08ex
\noindent%
\bfcour{\yblue Quantification of vector.}
Inasmuch as vector space is but a formal set of axioms \cite{halm},
let us take the premise that we have no primitive concepts beyond
those encompassed within the concept of a linear vector space (the
bare \lvs): no length/\!"!norm, no rotations/\!"!angles, operators, no inner
product, orthogonality, and no other non"=canonical terms or related
interpretations. Even the concepts of a basis and a model of \lvs\
can be dropped. While occasionally we invoke geometric arguments in
the \lvs"=mathematics, they can eventually be fully discarded.

Let us also pose the problem of `quantifying the vector', \ie,
creating a numerical quantity $\cal Q$ we attach to every vector
$\greekC \in V$. This is referred to as the \emph{quantitative
measure} of a vector or simply magnitude $\cal Q[\greekC]$.

The underlying reason to introduce such an object is that the vector
itself possesses no intrinsic quantitative
properties/\!"!characteristics---be they geometrical or, to use the
quantum"=theoretical lexicon, measurable, observable, or
interpretative \cite{aaronson}. Vector coordinates cannot serve as
such characteristics because the coordinate isomorphism $\ds
V\mapsto \bbR^n$ is not canonical. And not being the 1"=dimensional
objects, vectors are not computed like quantities;
multi"=dimensionality of \lvs\ is essential. What might be computed
are the 1"=dimensional magnitudes in their own right, those derived
from vectors. What properties, then, characterize $\cal Q$ as a
scalar add-on to the $V\!\!$"=space?

Once geometry is formalized in terms of an \lvs, its
multi"=dimensionality manifests in any vector being a
$\smallc({\+}\smallc)$"=sum of many other ones: $\msf c\bcdot\greekC
= \msf a\bcdot\greekA\,\,\+\,\, \msf b\bcdot\greekB$. This means that
the $\cal Q$"=theory must be built upon the coefficients $\{\msf a,
\msf b, \msf c\}$, for these numbers determine the variety of all
possible triangles in Euclidean geometry. The theory of magnitudes
$\cal Q[\msf c\bcdot\greekC]$,
$\cal Q\mbig[1][\msf a\bcdot\greekA\,\,\+\,\,
\msf b\bcdot\greekB\mbig[1]]$ should thus become a numerical theory
of these $(\msf a, \msf b; \msf c)$---the numeric representatives of
the abstract triangles ${(}\![5]{(}\msf a, \msf b; \msf c\,|\, \greekA,
\greekB; \greekC{)}\![5]{)}$. Hence, $\msf a,\msf b\in\bbR$ are
understood to be arbitrary throughout.

In other words, the $\cal Q$"=value as a number may not originate
from anything but coordinates despite their being arbitrary even for
a fixed vector. Consequently, we have no way of producing this new
number other than from the $\bbR$"=field numbers $\{\msf a,
\msf b,\ldots\}$, while duly accounting for the presence of the
$\greekA$"=symbols. This is what we mean by `quantifying a vector'.

Let us summarize the reasoning provided above. Triangles and vectors
do not reduce to numbers. This is why their computations should
differ from just the arithmetic of numbers per~se, followed by a
quantitative labeling thereof; \eg, a large $\msf b$"=coefficient in
$\msf a\bcdot\greekB \+ \msf b\bcdot\greekA$, a long/\!"!short length
$|\msf c|$ of vector $(\!-\msf c\,{\bcdot}\,\greekA)$, \etc. Given a
vector, its quantitative meaning/\!"!interpreting is not well"=defined
a~priori. Thus, in the linear theory, quantification and the
interpretable must be a theory in its own right.

\vskip0.2ex
\noindent
\itcour{\yblue Additivity,~\,\,counting,~\,\,and~\,\,semigroup.}~%
In analogy to the additive measure on a set \cite{kuzma}

\vbox to3.5ex{\vspace{-1.1ex}
\begin{equation*} \mu(A\cup B) = \mu(A) + \mu(B)\.,\quad \text{if\;\;} A\cap B
= \varnothing\.,
\end{equation*}}
\noindent
we begin with additivity as a natural property of $\cal Q$:

\vbox to4.ex{\vspace{-0.9ex}
\begin{equation}\label{arrow}
\cal Q\big[\msf a\bcdot\greekA \+ \msf b\bcdot\greekB\big] =
\cal Q[\msf a\bcdot\greekA\,] + \cal Q[\msf b\bcdot\greekB]\..
\end{equation}}
\noindent
The matter of the pair $\{\greekA, \greekB\}$ is discussed further
below. Indeed, quantification is usually identified with the notion
of counting; say, with such concepts as `area of, volume of,
number/\!"!size of', \etc\ \cite{br2}. On the other hand, the primary
property that formalizes the quantification in these examples is
additivity \eqref{arrow}. Note that the `counting of things' implies
nothing about taking one away from the other, \ie, subtraction. The
magnitude cannot be negative. We hence adopt the following
formalization.%
\hypertarget{D1}{}

\vskip0.1ex
\textbf{Definition~1.} \emph{The quantitative measure
\textup{(}\!quan\-tifi\-ca\-tion\textup{)} of vectors is a functional
superstructure to the $V\![3]$"=space, the semigroup
$\cal S\langle\cal Q\in\ds\bbR^{\!\sp}_0; +\rangle$, where $\cal Q =
\cal Q[\greekC]$ is a well"=defined nonconstant map $\ds V\mapsto
\bbR^{\!\sp}_0$ for all $\greekC\in V\!\!$.}

\vskip0.1ex
\hypertarget{R1}{\textbf{Remark~1.}} Unlike the abstract $\mu$"=measure,
$\cal Q$ is not just a function on \lvs\ to be found, or to be
defined, such as a norm $\|{\cdot}{\cdot}{\cdot}\|$ \cite{halm,
aaronson}. The $\cal Q[\cdot{\cdot}\cdot]$ is realized as an
a~priori unknown well"=defined map $\bo\aleph$ on the
coordinate"=representations of vectors:
\vspace{-1.7ex}
\begin{equation}\label{xi}
(\bbR^2 \Times V^2)
\mathrel{\Over[1.1ex]{\!\text{\footnotesize$\bo\aleph$}\,}
{\longmapsto}} \bbR^{\!\sp}_0{:}\qquad\bo\aleph(\msf a, \msf b; \greekA,
\greekB) = \cal Q[\msf a\bcdot\greekA \+ \msf b\bcdot\greekB]\..
\end{equation}

\vspace{-1.4ex}
\noindent
This induces a calculus on the $V\!\!$"=space, henceforth referred to
as $\cal Q$"=calculus, a calculus of the $\cal Q$"=measures.

\hypertarget{Inew}{}\smallskip\noindent%
\bfcour{\yblue Extending the $V\![4]$"=space algebra\!:~\.$\cal
Q$"=measure.}
As it turns out, the additive characterization of the
$\cal Q$"=measure \eqref{arrow} along with its other properties
admits a purely formal justification within a minimalist framework.
We now elaborate on this idea.

\smallskip
\noindent
\itcour{\yblue Vectorial and numerical additions.}
Let us forget additivity \eqref{arrow}. Although the
$\cal Q$"=objects are the $\ds\bbR^{\!\sp}_0\!$"=numbers, they \emph{are
not} elements of the field $\bbR$, the field of coefficients
$\{\msf a, \msf b, \ldots\}$ over which the $V\!\!$"=space has been
defined. This means that, besides the standard \lvs"=distributivity
$\msf a\bcdot\greekA \.\+\. \msf b\bcdot\greekA = (\msf a+\msf
b)\bcdot\greekA$, a new relation needs to be included to conform the
additions: the $\smallc({\+}\smallc)$ for vectors and new plus
$\smallc({+}\smallc)$ for $\cal Q$\,s. Otherwise, these $\cal Q$\,s
\emph{would have no relation} to \lvs. Once introduced, these pluses
must be compatible. What can one say about this?

Given the need to introduce a new $\smallc({+}\smallc)$"=addition,
let us find a statement about compatibility between the binary
structure~$\greekA\+\greekB$ and the binary numerical one
$\cal Q_\tgreek{a}+\cal Q_\tgreek{b}$ via a certain function
$\frak L[\cdot{\cdot}\cdot]$. That is, without loss of generality we
should declare the \lvs"=representatives
$\mbig[1]\{\frak L[\greekA\,]\FED\cal Q_\tgreek{a},\;
\frak L[\greekB]\FED\cal Q_\tgreek{b}\mbig[1]\}$ for the
$\cal Q$"=numbers being added~$(+)$ and an $\frak L$"=representation
of $(\greekA\+\greekB)$, \ie, the object $\frak L[\greekA \+
\greekB]$. Then we are interested in an algebraic implementation of
the formal axiomatics ($\scr A_{\raisebox{-0.1ex}{$\sp$}}$)

\vbox to5.7ex{
\begin{equation*}
\big\{\greekA, \greekB,\ldots; \+\big\}_\text{\!\!\lvs}\;\.
\Under[1.1ex]{{\sss\text{compatibility\;of\;%
(\!\!$\Over[1.4ex]{\scalebox{0.3}[0.7]{$\sss\bo\frown$}}
	{\vcenter{\hbox{$\sss+$}}},\![3] +\!\!)$}}}
{\Over[1.7ex]{\sss\text{introduce $\{\!\!\cal Q; +\!\}$}}
{\goto[4]\.\cdots\.\goto[4]}}\;\;\, \frak L[\greekA \+ \greekB] =
\Over[2.2ex]{\smash{\.\sss(?)}}{\scr A_{\raisebox{-0.1ex}{$\sp$}}}
\!(\cal Q_\tgreek{a}, \cal Q_\tgreek{b})
\end{equation*}}
\noindent
and in determining the vectors $\greekA, \greekB\in V$ for which
such an $\scr A_{\raisebox{-0.1ex}{$\sp$}}$ exists.

At this stage, the $\{\cal Q_\tgreek{a},
\cal Q_\tgreek{b}\}$"=numbers are not subtracted, nor are they
multiplied yet; only $\smallc({+}\smallc)$. Technically speaking, we
would have to adopt a new plus"=symbol for this single operation
in the algebra of $\cal Q$\,s.%
\hypertarget{L1}{}

\smallskip
\textbf{Lemma~1.} \emph{$\scr A_{\raisebox{-0.1ex}{$\sp$}}\!(\cal
Q_\tgreek{a}, \cal Q_\tgreek{b}) = \cal Q_\tgreek{a} + \cal
Q_\tgreek{b}$ for some $\greekAo,\greekBo\.{\in}\.V\!\!$.}

\smallskip
\noindent
\emph{Proof.} %
All that the operational minimalism of the semigroup $\cal S$ does
allow us to do with its elements is to add them up:
$\cal Q_\tgreek{a} + \cal Q_\tgreek{b}$. Hence, should the
$\scr A\!$"=compatibility $\frak L[\greekA \+ \greekB] =
\scr A_{\raisebox{-0.1ex}{$\sp$}}\!\!\mbig[1](\frak L[\greekA\,],
\frak L[\greekB]\mbig[1])$ exist for some vectors $\greekAo$ and
$\greekBo$, it necessarily takes the form of the linear expression
\begin{equation*}
\frak L[\greekAo \+ \greekBo] =
\frak L[\greekAo] + \frak L[\greekBo]\qquad
\big\lceil\,\,\Over{\frak L\.\goto\.\cal Q\.}{\goto[5]}\;\eqref{arrow}\,
\big\rceil\..
\tag*{$\blacksquare$}
\end{equation*}

\vspace{-0.5ex}
\hypertarget{R2}{\textbf{Remark~2.}} We do not assume, a~priori, that
$\{\greekAo,\greekBo\}$ are arbitrary. Hence, $\frak L$ should not
be thought of as a homomorphism $(\+)
\mathrel{\Over[1.1ex]{\!\frak L\,}{\mapsto}} (+)$, \ie, as a linear
functional~\cite{halm}
\vbox to5ex{
\vspace{-0.3ex}
\begin{equation}\label{Lf}
\scr L[\greekA \+ \greekB] = \scr L[\greekA\,] + \scr L[\greekB]\qquad
\forall\, \greekA, \greekB\in V.
\end{equation}}
The quantitative theory of \lvs"=triangles is determined by numbers
$\{\msf a, \msf b\}$ rather than by non"=numerical $\{\greekA,
\greekB\}$"=objects in \eqref{Lf}. The latter are merely the formal
symbols without any meaning like length. This, together with
\hyperlink{L1}{\black Lemma~\brown\textbf{1}}, allows us to view the
additivity \eqref{arrow} almost as a derived property, not as an
ad~hoc postulate.

\smallskip
In what cases does the $\cal Q$"=measure consistently exist? Are
vectors $\{\greekAo,\greekBo\}$ arbitrary, (non)parallel? And in
general, is (non)parallelism $\greekA \nparallel \greekB$ an
(un)necessary condition for \eqref{arrow}? By the (geometric)
parallelism $\greekA\parallel\greekB$ is meant the formalization
$\greekB = \msf k\bcdot\greekA$ with some~$\msf k \in
\bbR\,\vcenter{\Smaller{\bo{\backslash}}}\, \{\msf0\}$. Let us call
$(\greekAo,\greekBo)$ the additive pair or, in short, the
$\scr A\!$"=pair. 

\smallskip\noindent
\itcour{\yblue Multiplication and endomorphisms.}
What can one say about the second operation on $V$? That is, how
does $\cal Q$ change if vector is multiplied by a constant:

\vbox to4.3ex{
\vspace{-0.8ex}
\begin{equation}\label{mult}
\smash{\cal Q[\greekA\,]\quad \longmapsto\quad
\cal Q[\msf c\bcdot\greekA\,] = {}\bfcour?}\..
\end{equation}}

Were we in elementary geometry, we could write
$\cal Q[\msf c\bcdot\greekA\,] = |\msf c|\times\cal Q[\greekA\,]$ for
length, $\msf c^2\times\cal Q[\greekA\,]$ for area, and
$|\msf c|^3\times\cal Q[\greekA\,]$ for volume""---""candidates for
the additive magnitude. However no such terms are available yet for
us. Besides, there is no reason why an additive quantity on \lvs\
should also be multiplicative: $\cal Q[\msf c\bcdot\greekA\,] =
\hbox{\small$\mco{const}$}\times\cal Q[\greekA\,]$. It does not
follow from anywhere. Incidentally, the standard length \cite{halm}
is clearly not additive, so this first candidate is inherently
unsuitable.

Much like the interplay between pluses, the consistent theory
necessitates a compatibility between the internal symmetries
$V\mathrel{\Over[1.4ex]
{\!\smallc[3](\msf c\bcdot\!\!\smallc[4])}{\mapsto}}V$ of the
$V\!\!$"=algebra and the functional equivalence $\mbig[1]\{\greekA
\mathrel{\Over[1.2ex]{\sss\cal Q}\sim} \greekB \hhence
\cal Q[\greekA\,] = \cal Q[\greekB]\mbig[1]\}$ being induced by the
external map $V \mathrel{\Over[1.2ex]{\!\sss\cal Q\,\,}{\mapsto}}
\ds\bbR^{\!\sp}_0$. The necessary way to ensure that compatibility is
to require preservation of the $\cal Q$"=equivalence---a
$\cal Q$"=class gets $\smallc(\msf c{\bcdot}\!\smallc)$"=mapped to a
class:

\vbox to4.5ex{\vspace{-0.5ex}
\begin{equation}\label{well}
\big(\cal Q[\greekA\,] = \cal Q[\greekB]\big) \quad\hence\quad
\big(\cal Q[\msf c\bcdot\greekA\,] =
\cal Q[\msf c\bcdot\greekB]\big)\quad\forall\,(\msf c,\greekA).
\end{equation}}
\noindent
This invariance of the $\cal Q$"=fibres may be regarded as a
requirement for the $\cal Q$"=measure to be meaningful.%
\hypertarget{L2}{}

\textbf{\!Lemma~2.}~\emph{The functional relation holds}:

\vbox to4.5ex{
\vspace{-0.5ex}
\begin{equation}\label{Cc} \cal Q[\msf c\bcdot\greekA\,] =
\mco C(\cal Q[\greekA\,])\qquad \forall\,(\msf c, \greekA\,)%
\in (\bbR\Times V)\,.
\end{equation}}%
\noindent
\emph{Proof.} Consider operator automorphisms
$\greekA\to\skew1\hat{\msf c}\,\greekA\FED\msf c\bcdot\greekA$ of the
$V\!\!$"=space as a commutative group \cite{halm}. Each of them
induces, due to \eqref{well}, a well"=defined map on
$\ds\bbR^{\!\sp}_0$:
\vbox to5ex{
\vspace{-0.5ex}
\begin{equation*}
\big(V\ni\greekA\to\skew1\hat{\msf c}\,\greekA\in V\big)\quad
\hence\quad \big(\bbR^{\!\sp}_0\ni\cal Q[\greekA\,]\to
\cal Q[\skew1\hat{\msf c}\,\greekA\,]\in\bbR^{\!\sp}_0\big)\,.
\end{equation*}}
\noindent
However such relation between number sets---$\{\cal Q[\greekA\,]\}$
and $\{\cal Q[\skew1\hat{\msf c}\,\greekA\,]\}$---is a
numerical~function. Hence it follows \eqref{Cc}. \hfill
{$\blacksquare$}

The function~$\mco C$ depends on $\msf c$, of course. See also the
topological caveats in \hyperlink{R3}{\black Remark~\brown\textbf{3}}
further below.

\smallskip
\noindent
\itcour{\yblue Bases and $\boB\![5]$"=invariance of $\cal Q$.}
Bases form a set $\scr B\ni(\greekA, \greekB)$ in its own right,
where $\greekB\ne\msf k\bcdot\greekA$. In turn, vectors, magnitudes,
and the $\scr B$"=set are independent, in a sense, as the
set-theoretical structures $\{\text{\lvs}, \scr B\!\!, \cal S\}$
impose no constraints on each other. Therefore, we should adopt both
the basis and the $\cal Q$"=numerical freedom:\hypertarget{boB}{}

\vskip1ex
\noindent
\hfil$(\boB\!)\;\;$%
\parbox[t]{6.5cm}{Bases $\scr B$ and the $\cal Q$"=measures of vectors
are independent of one another.}

\vskip1ex
\noindent
One aspect of this independence is the standard invariance
\cite{halm, kuzma}: the $\cal Q$"=value for a given vector should not
change under a change of basis $(\greekA, \greekB)\goto (\greekA',
\greekB')\in\scr B$.

Let us take a closer look at the choice of a $\scr B$"=element. We
could, for example, ab~initio declare the pair $\{\greekA,
\greekB\}$ to be the basis $(\greekB, \greekA\,)$. Such a re"=labeling
`first $\rightleftarrows$ second' of the ordered symbols $(\greekA,
\greekB) \rightleftarrows (\greekB, \greekA\,)$ must not change the
$\cal Q$"=value because neither of these orders is pre"=defined.
Indeed, which of the two bases should be associated with the
triangle ${(}\![5]{(}\msf a,\msf b;\msf c\,|\ldots {)}\![5]{)}$, \ie, how
should we label the triangle with vector"=symbols $\{\greekA,
\greekB, \greekC\}$? Why $(\greekA, \greekB)$ and not $(\greekB,
\greekA\,)\in\scr B$, given that $\scr B$ is only a set? Expressed
another way, when geometry of triangles is reformulated in the
abstract \lvs"=terminology, the pair $\{\msf a\bcdot\greekB,
\msf b\bcdot\greekA\}$ might represent triangles equally well as
does the seemingly obvious pair $\{\msf a\bcdot\greekA,
\msf b\bcdot\greekB\}$. We will call the mutual
$\{\cal Q,\scr B\}$"=independence the $\boB\!\!$"=\emph{invariance}
(`\textbf{b}ases \tplus\ $\cal Q$').

Formally, the permutation $(\greekA, \greekB) \rightleftarrows
(\greekB,\greekA\,)$ is the natural automorphism on the $\scr B$"=set.
Thus, we must impose an equivalence relation $(\greekA, \greekB)
\sim (\greekB,\greekA\,)$ on $\scr B$, which entails the permutation
equality

\vbox to4.5ex{\vspace{-0.5ex}
\begin{equation*}
\cal Q\big[\msf a\bcdot\greekA \+ \msf b\bcdot\greekB\big] =
\cal Q\big[\msf a\bcdot\greekB \+ \msf b\bcdot\greekA\big]\qquad
\forall\, \msf a, \msf b\.,
\end{equation*}}
\noindent
even though $\msf a\bcdot\greekA\,\,\+\,\, \msf b\bcdot\greekB \ne
\msf a\bcdot\greekB\,\,\+$ $\msf b\bcdot\greekA\,$.%

\vskip0.3ex\noindent%
\hypertarget{I3}{}%
\bfcour{\yblue Definition $\goto$ Theorem.}
As noted above, these principles are sufficient to derive the
Pythagorean theorem itself and the geometric framework surrounding
it \cite{maor, givental}. This minimalism also provides the core of
a rigorous justification for the (`quantum') complex field $\bbC$.

Given what has been said about triangles and \lvs, together with
\hyperlink{L1}{\black Lemmas~\brown\textbf{1}}--%
\hyperlink{L2}{\brown\textbf{2}} and \hyperlink{R2}{\black
Remark~\brown\textbf{2}}, let us formulate the structural properties
determining the $\cal Q$"=object. Although they are derivable (see
\hyperlink{T2}{\black Theorem~\brown\textbf{2}} further below), we
temporarily organize them, for convenience, into a set of technical
axioms. Moreover, the assumption that $(\greekA, \greekB)$ be a
basis can be dropped here.%
\hypertarget{D2}{}

\smallskip
\textbf{Definition~2.} \emph{Axioms of the $\boB\![4]$"=invariant
quantitative measure $\cal Q$ on coordinate"=representations of
vectors}:
\newcounter{myQ}%
\newcounter{myA}%
\begin{flalign}
\Over[2.2ex]{\smash{\sss(?)\!}}{\cal Q} [\msf c\bcdot\greekC]\,[1.5]
&= \smash{\Over[2ex]{\smash{\sss(?)\,}}{\mco C_{\msf c}}}
(\cal Q[\greekC])&\makebox[0ex][r]{$\forall\,(\msf c,\greekC)$}&%
\makebox[0ex][l]{${}\in (\bbR\,\,{\Times}\,\,V)\,,$}
\label{3}\\
\cal Q\big[\msf a\bcdot\greekA \+ \msf b\bcdot\greekB\big] & = \cal
Q{\,[1.5]}[\msf a\bcdot\greekA\,] + \makebox[0ex][l]{$\cal
Q{\,[1.5]}[\msf b\bcdot\greekB]$}& \forall\,\msf a, \msf
b&\in\bbR\.,\qquad\.
\label{4}\setcounter{myQ}{\value{equation}}
\\
\cal Q\big[\msf a\bcdot\greekA \+ \msf b\bcdot\greekB\big] & = \cal
Q\big[\msf a\bcdot\greekB \+ \msf b\bcdot\greekA\big] &
\forall\,\msf a, \msf b&\in\bbR\..
\label{4*}\setcounter{myA}{\value{equation}}
\end{flalign}

\vspace{-1ex}
\noindent
\emph{Vectors} $\greekA,\greekB$ \emph{are to be determined by
requiring the consistent existence of\ \,$\cal Q$ and
are~called~the~admissible,~$\bo\aleph$"=pair.}

\smallskip
\noindent
\itcour{\yblue$\mco Q$"=representation.}
Reframing the task, what is being sought is not merely functions
$\cal Q$ and $\mco C_{\msf c}$ but also the existence domain for the
$\{\greekA, \greekB\}$"=arguments of the formal 4"=argument
$\bo\aleph$"=function \eqref{xi} under arbitrary $\msf a$ and
$\msf b$. Naturally, the $\greekA$"=symbols of the vector abstracta
will, one way or another, have to disappear in the ultimate result:
\begin{equation}\label{star}
\cal Q\big[\msf a\bcdot\greekA \+ \msf b\bcdot\greekB\big]
\quad\smash{\Over[1.3ex]{\ds\bfcour?}{\goto[3]}}\quad
\big\lceil\,\smash{
\Over[2.4ex]{\sss(?)}{\mco Q}(\msf a,\msf b)\.,\quad
\{\Over[2.4ex]{\sss(?)}{\greekA, \greekB}\}}\, \big\rceil\..
\end{equation}
The $\mco Q$ here is a numeric representative of
$\cal Q[{\cdot}{\cdot}{\cdot}]$. This raises the question of the
procedure \eqref{star}, \ie, of whether such representations do
exist and, if so, what is the mechanism for obtaining the purely
numerical form $\mco Q({\cdot}{\cdot}{\cdot})$ free of symbols
$\{\greekA,\+,\bcdot\}$, a `$\mco Q$"=theory' of coordinates.%
\hypertarget{R3}{}

\smallskip
\textbf{Remark~3}~(topology). In this note, we are interested only in
the algebraic aspect of the problem. We also adopt the standard
topological conventions \cite{kuzma} and use them wherever
necessary. Accordingly, we require $V$ to be equipped with a
topology compatible with the linear structure, so that
continuity/\!"!analyticity is defined for the functions that will arise
or be implicitly assumed in what follows: $\cal Q[\greekA\,]$,
$\mco C(\cal Q[\greekA\,])$, $\mco Q({\cdot}{\cdot}{\cdot})$, \etc.
Since $V$ is finite"=dimensional, it carries a unique Hausdorff
vector topology. This ensures the existence of these functions as
continuous ones, provided the algebraic structure is consistent.

\smallskip
\hypertarget{T1}{\textbf{Theorem 1.}} \emph{Let $V$ be a
$2$"=dimensional vector"=space over $\bbR$. Then the $\cal Q$"=measure
of vector, as defined above, has the following
$\mco Q$"=representation}

\vspace{-0.5ex}
\vbox to5ex{
\begin{equation}\label{Q}
\cal Q\big[\msf a\bcdot\greekA \+ \msf b\bcdot\greekB\big] =
\itcour{const}\times(\msf a^2+\msf b^2)\..
\end{equation}
}
\noindent
\emph{Here, the admissible $\{\greekA, \greekB\}$"=pairs must form a
basis and are related to each other through an orthogonal
transformation and determine the right"=angled triangles
$\bo[\hskip-1.6pt\bo[\msf a, \msf b;
\msf c\,|\,\greekA,\greekB;\greekC\bo]\hskip-1.6pt\bo]$. For such
triangles, the relation $\msf a^2+\msf b^2 = \msf c^2$ holds as a
gauge"=invariant identity.}

The concepts of `orthogonal, right"=angled, and gauge` will naturally
arise and be defined in the course of the proof. It will take, in
effect, the form of a derivation and reveal the basis changes to be
partial. The function $\itcour{const} = \itcour{const}[\greekA,
\greekB]$ is also to be determined.%

\smallskip\noindent%
\hypertarget{I4}{}%
\bfcour{\yblue\!How~\,\,does~\,Pythagorean~\,square~\,come~\,\,%
into~\,\,being\!?} %
Let us proceed to the proof of \hyperlink{T1}{\black
Theorem~\brown\textbf{1}}. \hypertarget{I}{}%

\vskip0.3ex
\noindent
\itcour{\yblue\textup{\![3](\!I\!)}~~Functionality of $\cal Q$.}
Suppose there exists a nonconstant continuous function
$\cal Q[\greekC]$ for all $\greekC\in V$. Take an arbitrary
(initial) pair $(\greekAo, \greekBo)$ and declare it to be an
$\bo\aleph$"=pair, \ie, suppose the expressions
$\cal Q[{\cdot}{\cdot}{\cdot}]$ satisfy \eqref{3}--\eqref{4*}. Then,
dropping the subscript {\footnotesize{$\msf o$}}, the distributive
and associative laws from axioms of \lvs\ \cite{halm} entail the
following identities:
\begin{align*}
\msf c\bcdot(\msf a\bcdot\greekA \+ \msf b\bcdot\greekB)\,&
\EqDown{}{\footnotesize\qquad\ceils{$\msf c\bcdot(\msf a\bcdot\greekA\,)
= \msf{c\,a}\bcdot\greekA\,$}}\,\msf c\bcdot(\msf a\bcdot\greekA\,) \+
\msf c\bcdot(\msf b\bcdot\greekB)\qquad
\forall\,\msf a, \msf b, \msf c\\
\cal Q\big[\msf c\bcdot(\msf a\bcdot\greekA \+
\msf b\bcdot\greekB)\big]&
\EqDown{}{\footnotesize\qquad\ceils{$\eqref{4},\;\; \msf{c\,a}\bcdot\greekA
= \msf c\bcdot(\msf a\bcdot\greekA\,)$}}
\cal Q\big[\msf{c\,a}\bcdot\greekA \+ \msf{c\,b}\bcdot\greekB\big]\\
\cal Q\big[\msf c\bcdot(\msf a\bcdot\greekA \+
\msf b\bcdot\greekB)\big]&
\EqDown{\footnotesize\ceils{\eqref{3}}\qquad}
{\footnotesize\qquad\ceils{$\eqref{3},\;\;
\mco C_{\msf c}\FED\mco C$}}
\cal Q\big[\msf c\bcdot(\msf a\bcdot\greekA\,)\big] +
\cal Q\big[\msf c\bcdot(\msf b\bcdot\greekB)\big]\\
\mco C\big(\cal Q[\msf a\bcdot\greekA \+ \msf b\bcdot\greekB]\big)&
\EqDown{\footnotesize\ceils{\eqref{4}}\qquad}{}
\mco C\big(\cal Q[\msf a\bcdot\greekA\,]\big) +
\mco C\big(\cal Q[\msf b\bcdot\greekB]\big)\\
\mco C\big(\cal Q[\msf a\bcdot\greekA\,] +
\cal Q[\msf b\bcdot\greekB]\big)&\!\LongEqDown[4ex]
{\small$\mbig[0]\lceil$\footnotesize\!{$\cal Q[\msf a\bcdot\greekA\,]
\FED x\ne\const$\small$\mbig[0]\rceil$}\quad}
{\!\!\quad\small$\mbig[0]\lceil$\footnotesize\!$\cal Q[\msf b\bcdot\greekB]
\FED y\ne\const$,\;\;\hyperlink{R3}{\black Remark~\brown\textbf{3}}%
\small$\mbig[0]\rceil$}
\mco C\big(\cal Q[\msf a\bcdot\greekA\,]\big) +
\mco C\big(\cal Q[\msf b\bcdot\greekB]\big)\quad\forall\,\msf a,\msf b
\end{align*}
\vbox{
\vskip-4ex
\begin{equation*}
\mco C(x+y) \mathbin{\Over[1.7ex]{\forall\, x\!,y}{=}} \mco C(x) +
\mco C(y)\quad\Over[1.7ex]{\text{\cite{hamel, kuzma}\.}}
{=\![3]\Longrightarrow}\quad\mco C(x) = \hbox{\small$\mco{const}$}\times
x\..
\end{equation*}}

\vspace{-1.5ex}
\noindent
Therefore, the semigroup $\cal S$ should be equipped with a
numerical $\smallc({\times}\smallc)$"=multiplication operation:

\vspace{-0.5ex}
\vbox to 5ex{
\begin{equation*}
\cal Q[\msf c\bcdot\greekA\,] =
\hbox{\small$\mco{const}$}(\msf c)\times\cal Q[\greekA\,],\qquad
\cal Q[\msf c\bcdot\greekB] = \cdots \quad\forall\,\msf c\..
\end{equation*}}

\noindent
Furthermore, these expressions and \eqref{3} imply
$\mco C_{\msf c}(\cal Q)=
\hbox{\small$\mco{const}$}(\msf c)\times\cal Q$ for all $\cal Q$\,s.
The rule \eqref{mult}--\eqref{3} is hence the
$\ds\bbR_{\smash{\sss\!\times}}^{\!\sp}$"=multiplication indeed:

\vbox to5.5ex{
\begin{equation}\label{3+}
\cal Q[\msf c\bcdot\greekC] =
\hbox{\small$\mco{const}$}(\msf c)\times\cal Q[\greekC]
\qquad\forall\,(\msf c,\greekC)\in (\bbR\Times V)\..
\tag{$\green\textsf{\themyQ}'$}
\end{equation}}

\noindent
This also yields the $(\PM+)$"=involution (see \cite[p.~13]{br2}):
\begin{gather}
\cal Q[\greekC] = \cal Q[-(-\greekC)] = \big|\eqref{3+}\big| =
\hbox{\small$\mco{const}$}^2(-\msf1) \times \cal Q[\greekC] \quad\hence
\notag\\[1ex]
\hbox{\small$\mco{const}$}(\pm\msf1) = +1\quad\hence\quad \cal Q[\pm\greekC] =
\cal Q[\greekC]\label{pm}
\end{gather}
because there are no negatives in semigroup.

If $\greekA\parallel\greekB$ then, putting
$\greekB=\msf k\bcdot\greekA$, we obtain an answer to the above
question on (non)parallelism. Indeed, \big|\eqref{4}\big| \hence\
$\cal Q[\greekA+\msf k\bcdot\greekA\,] = \cal Q[\greekA\,] +
\cal Q[\msf k\bcdot\greekA\,]$ \hence\ $\big|\eqref{3+}\big|$ \hence\
$\hbox{\small$\mco{const}$}(\msf1+\msf k) =
1+\hbox{\small$\mco{const}$}(\msf k)$ \hence\
$\hbox{\small$\mco{const}$}(\msf k)=\msf k\;\;\forall\, \msf k\in\bbR$
\cite{kuzma, hamel}; con\-tra\-dic\-tion, because $\msf k$ may be
negative. Accordingly, for a nontrivial
$\cal Q[{\cdot}{\cdot}{\cdot}]$, the $\scr A\!$"=pair
$(\greekAo,\greekBo)$ must form a basis.
We are thus interested in, and are looking for, the form\\
\vbox to6ex{
\begin{equation}\label{Qa+}
\cal Q[\msf a\bcdot\greekC] =
\Over[2ex]{\smash{\sss(?)\!}}{\mco Q}(\msf a)\times
\cal Q[\greekC]\qquad\forall\, \msf a, \greekC\.,
\tag{$\green\textsf{\themyQ}''$}
\end{equation}}

\noindent
wherein $\mco Q({\cdot}{\cdot}{\cdot})$ is a new function to be
found, the very numeric representative of
$\cal Q[\msf a\bcdot\greekC]$ in \eqref{star}. So, the linear
$V\!\!$"=structure entails separability of variables in
$\bo\aleph(\msf a,\msf b;\greekA, \greekB)$.

The $\smallc({\PM+}\smallc)$"=symmetry \eqref{pm} says that $\mco Q$
is an even function: $\mco Q(\msf a) = \varkappa_0 +
\varkappa_2\,\msf a^2 + \varkappa_4\,\msf a^4 + \cdots$. The
analyticity requirement (\hyperlink{R3}{\black
Remark~\brown\textbf{3}}) implies that there is no loss of generality
here; every analytic function is a series. The multiplicative
property \eqref{Qa+} becomes the homogeneity condition
$\mco Q(\msf{c\,a}) = \mco Q(\msf c)\times \mco Q(\msf a)$ and,
thereby, `kills' all terms in this series, save one: $\mco Q(\msf a)
= \msf a^{2\textsf{p}}$ \cite{kuzma}.

Consider now axiom \eqref{4*}. Setting $\msf b = \msf0$ and $\msf a
= \msf1$ there, one obtains a property the $\bo\aleph$"=basis must
satisfy:

\vbox to 5ex{
\vspace{-1.2ex}
\begin{equation}\label{QaQb}
\cal Q[\greekB] = \cal Q[\greekA\,]\.,\qquad
\cal Q[\msf a\bcdot\greekA\,]=\msf a^{2\textsf{p}}\times
\cal Q[\greekA\,]\..
\notag
\end{equation}}
\noindent
Therefore, in the framework of the topological conventions above, we
arrive at the intermediate result
\begin{equation}\label{Q2p}
\cal Q\big[\msf a\bcdot\greekA \+ \msf b\bcdot\greekB\big] =
(\msf a^{2\textsf{p}} + \msf b^{2\textsf{p}})\times \cal Q[\greekA\,]\quad
\forall\,\msf a,\msf b\in\bbR\.,
\end{equation}
where $\textsf{p} = 1,2,\ldots$ is an as yet unknown integer independent of
the basis $(\greekA, \greekB)$.

Up to this point we have dealt with a fixed $\scr A\!$"=pair
$(\greekA, \greekB)$. Consider next the
\hyperlink{boB}{$\brown\boB\!\!$}"=invariance.
\hypertarget{II}{}%

\smallskip\noindent
\itcour{\yblue\textup{\![3](\!I\!I\!)} Whence the two and gauge?}
Given a vector,
\vspace{-3.2ex}
\begin{align}
\msf a\bcdot\greekA \+ \msf b\bcdot\greekB\.& = \msf
a'\bcdot\greekA\,' \+ \msf b'\bcdot\greekB'\quad\;\;
\forall\,\msf a,\msf b\.,\notag\\[-1.5ex]
\intertext{and any other $\bo\aleph$"=bases $(\greekA', \greekB')
\in \scr B$, should they exist, have equal standing in $\scr B$.
They yield the same result:}%
\cal Q\big[\msf a\bcdot\greekA \+ \msf b\bcdot\greekB\big]\![1.5]&
\,\EqDown{\footnotesize $\text{\ceils{\eqref{QaQb}}}$\qquad}
{\footnotesize\qquad $\text{\ceils{\eqref{QaQb}}}$}
\cal Q\big[\msf a'\bcdot\greekA\,' \+ \msf b'\bcdot\greekB'\big]
\notag\\[0.5ex]
(\msf a^{2\textsf{p}} + \msf b^{2\textsf{p}}){\,\times\,}\cal
Q[\greekA\,] & = \big((\msf a')^{2\textsf{p}} + (\msf
b')^{2\textsf{p}}\big){\,\times\,}\cal Q[\greekA'] \hbox
to0ex{\.,}\label{p2}
\end{align}
where pairs $(\msf a,\msf b)$ and $(\msf a',\msf b')$ are related
through a linear $\big(\smallmatrix
a&b\\c&d\endsmallmatrix\big)$"=transformation $\{\msf a' =
a\,\msf a+b\,\msf b$,\; $\msf b' = c\,\msf a+d\,\msf b\}$\,.

The \hyperlink{boB}{$\brown\boB\!\!$}"=principle states that the basis
changes $(\msf a,\msf b)\goto (\msf a',\msf b')$ must not depend on
the superstructure $\cal Q$, \ie, on the values
\{$\cal Q[\greekA\,]$, $\cal Q[\greekA']$, \ldots\}. This forces us to
put in \eqref{p2} the equalities
\begin{equation}\label{Qa'}
\cal Q[\greekA\,] = \cal Q[\greekA'] =
\cal Q[\greekA''] = \cdots=\cal Q[\greekB] =\cdots\;\;(\text{gauge})\..
\end{equation}

\textbf{Remark~4.} We see that the $\cal Q$"=measures of all admissible
basic $\greekA$"=vectors must coincide with each other. Although the
genesis of this formal property is not obviously inferable, it is
itself a well-known fact addressed in standard texts. It is referred
to as a normalization of bases \cite{halm}. A natural analogy to
such a gauge freedom in the \ceils{\lvs\ \tplus\ $\cal Q$}"=theory is
the (multiplicative) change of physical units. Incidentally, in
quantum theory \cite{aaronson}, one puts a stricter normalizing
convention~$\cal Q[|\bo\alpha\rangle] = 1$.

Clearly, equalities \eqref{Qa'} cannot be valid for all $\greekA\in
V$. Then for which ones? The cancelation
$(\cdots)\times\cal Q[\greekA\,] \hspace{0.25em}
\smash{\raise0.1ex\hbox{$\llap{\rotatebox{75}{$\bigg|$}}$}\![3]} =
(\cdots)\times\!\cal Q[\greekA']\hspace{0.1em}\,
\smash{\raise0.1ex\hbox{$\llap{\rotatebox{75}{$\bigg|$}}$}}$ turns
\eqref{p2} into a restriction on $\{a,b,c,d\}$:
\begin{equation*}
\msf a^{2\textsf{p}} + \msf b^{2\textsf{p}} =
(a\,\msf a+b\,\msf b)^{2\textsf{p}} + (c\,\,\msf a+d\,\msf b)^{2\textsf{p}}\qquad
\forall\,\msf a, \msf b\in\bbR\..
\end{equation*}
By expanding and collecting the coefficients in front of binomial
terms like $\msf a^{2\textsf{p}\sm k}\msf b^k$, one splits the
restriction:

\vspace{-0ex}
\vbox to12ex{
\begin{multline*}
\big\{a^{2\textsf{p}\sm1}\,b^1 + c^{2\textsf{p}\sm1}\,d^1 = 0\.,\quad
\wavy[4]{a^{2\textsf{p}\sm2}\,b^2 + c^{2\textsf{p}\sm2}\,d^2 = 0}\.,\\[1ex]
\ldots\ldots\.;\quad a^{2\textsf{p}} + c^{2\textsf{p}} = 1,\quad b^{2\textsf{p}} +
d^{2\textsf{p}} = 1\big\}\..
\end{multline*}}

\noindent
Solution to these equations is almost obvious:
\begin{align*}
\textsf{p} = 1{:}&\;\big\{a\,b+c\,d = 0,\quad a^2+c^2=1,\quad
b^2+d^2=1\big\}\.,
\\[1ex]
\textsf{p} > 1{:}&\;\big\{\wavy{(a^{\textsf{p}\sm1}b)^2 +
(c^{\textsf{p}\sm1}d)^2 = 0},\;\ldots \big\}\hence \{a\,b = 0 =
c\,d\}\..
\end{align*}
The wavy"=terms imply, up to $\pm1$, only the identical transformation
$\big(\smallmatrix a&b\\c&d\endsmallmatrix\big) = \big(\smallmatrix
1&0\\0&1\endsmallmatrix\big)$ and the trivial permutation
$\big(\smallmatrix a&b\\c&d\endsmallmatrix\big) = \big(\smallmatrix
0&1\\1&0\endsmallmatrix\big)$. %
The corresponding $\cal Q$"=calculi are degenerate because they
represent the separated
$\cal Q^{{\sss(}\textsf{p}{\sss)}}_{\tgreek{a}\tgreek{b}}$"=`theories'
for each $(\greekA, \greekB)$"=element of the $\scr B$"=set. The case
$\textsf{p}>1$ is hence trivial in the sense that it admits no true
changes of bases. This contradicts the
\hyperlink{boB}{$\brown\boB\!\!$}"=invariance, and we reject this
choice. Only the case $\textsf{p} = 1$ remains nontrivial. It preserves
the sum of squares $(\msf a^2+\msf b^2)$ and thereby yields the
well"=known objects: $\big(\smallmatrix
a&b\\c&d\endsmallmatrix\big){\; =\;}
\big(\!\smallmatrix%
\hfill\cos\vartheta&\sm\!\!\sin\vartheta\hfill\\
\hfill\sin\vartheta&\phantom{\sm}\!\!\cos\vartheta\hfill
\endsmallmatrix\!\big){\;\goto\;} \bb O_2$. This is where the
concept of rotation emerges in geometry; in its \lvs"=reformulation,
to be precise. The notion of an angle does not appear here yet, nor
is it needed.\hypertarget{III}{}%

\smallskip
\noindent
\itcour{\yblue\textup{\![3](\!I\!I\!I\!)}~Orthogonal bases and right
triangle.}
Thus, declaring the arbitrary initial $(\greekAo, \greekBo)$ to be
an $\bo\aleph$"=pair is consistent. Moreover, the pair belongs to a
family $\scr B_{\raisebox{0.1ex}{$\!\sss\msf o$}}\subsetneq\scr B$,
for which only continuous $\vartheta$"=rotations and discrete
reflections/\!"!permutations are allowed as the basis changes. We adopt,
by convention, the term `orthogonal' for these special bases. In
analogy with the above, no concept of orthogonality on $V$ is
required at this stage. Clearly, these bases form a class of
orthogonal equivalence, and one of its representative may be, as we
have seen, any element $(\greekA, \greekB)\in\scr B$. Returning to
geometry, we call the resulting triangles the right"=angled ones
$\hbox{$\bo{\mbig[0]{[}}\hskip-2pt\bo{\mbig[0]{[}}$} \msf a,
\msf b;...|\,\, \greekA\,,\greekB\,;...
\hbox{$\bo{\mbig[0]{]}}\hskip-2pt\bo{\mbig[0]{]}}$}$ and introduce
the notion of perpendicular $\greekA \mathrel{\bot} \greekB$.

The previous reasoning has not addressed the right hand side of
equality $\msf a\bcdot\greekA\,\,\+\,\, \msf b\bcdot\greekB =
\msf c\bcdot\greekC$. \hbox{On the other hand}

\vbox to4.5ex{
\vspace{-0.4ex}
\begin{equation}\label{delta}
\msf a \bcdot\greekA \+ \msf b\bcdot\greekB =
\msf c \bcdot\greekC \+ \msf0 \bcdot\greekD =
\msf c'\bcdot\greekC'\+ \msf0 \bcdot\greekD'= \cdots\.,
\notag
\end{equation}}
\noindent
where $\{\greekC, \greekD\}$\,s are also bases. The question then
arises: Can $\greekC$ be an element of an orthogonal basis and,
accordingly, satisfy $\cal Q[\greekC] = \cal Q[\greekA\,]$, followed
by $\cal Q[\msf c\bcdot\greekC]=\msf c^2 \times \cal Q[\greekC]$?
The answer is of course yes because any vector, including
$\{\greekC, \greekC', \ldots\}$, can be rotated up to a multiplicative
constant into any vector $\greekA$. From \eqref{Qa'} it follows that
we can always equate $\cal Q[\greekC] = \cal Q[\greekA\,]$ for one of
the $\greekC$\,s. As in \eqref{QaQb}, this amounts to assigning a
$\cal Q$"=value to~$\greekC$.

By invoking the equalities $\cal Q[\greekA\,] = \cal Q[\greekB] =
\cal Q[\greekC]$, we are now in a position to implement the scheme
\eqref{star}. The disappearance of the \lvs"=abstracta $\{\greekA,
\+, \bcdot\}$ takes the form
\vspace{-3.5ex}
\begin{align*}
\cal Q\big[\msf a\bcdot\greekA \+ \msf b\bcdot\greekB\big]&
\,\EqDown[3.5ex]{}{} \cal Q[\msf c\bcdot\greekC]
&&\hskip5pt\lower0.2ex\hbox{$\Big|$}\begin{array}{@{}l@{}}
\text{\smaller\;right-angled triangles}\\
\text{\smaller\;$\bo[\hskip-1.6pt\bo[\msf a,
\msf b; \msf c\,|\, \greekA, \greekB;\greekC\bo]\hskip-1.6pt\bo]$}
\end{array}\\[-1.9ex]
\;\;\;\;(\msf a^2 + \msf b^2)\times\,\,\cal Q[\greekA\,]\hspace{0.07em}
\smash{\llap{\rotatebox{162}{$\rule{21pt}{0.5pt}$}}}
& = \msf c^2\times\,
\cal Q[\greekC]\hspace{0.07em}
\smash{\llap{\rotatebox{162}{$\rule{21pt}{0.5pt}$}}}
&&\hskip5pt\Big|\begin{array}{@{}l@{}}
\text{\smaller\;full vector form of}\\
\text{\smaller\;Pythagorean theorem}
\end{array}\\[-2.2ex]
\intertext{and yields the gauge"=invariant numeric relation}\\[-5.4ex]
\mco Q{:}\qquad\msf a^2+\msf b^2& = \msf c^2\,.
&&\hskip5pt\Big|\smash{\begin{array}{@{}l@{}}
\text{\smaller\;Pythagoras; the numeric}\\
\text{\smaller\;representation of
$\cal Q[{\cdot}{\cdot}{\cdot}]$}
\end{array}}
\end{align*}

\vskip-0.6ex
\noindent
Observe that $\smallc({+}\smallc)$ in this gauge"=invariant
`$\mco Q$"=theory' \eqref{star} is the semigroup plus. The overall
scale $\cal Q[\greekAo] = \itcour{const}[\greekAo, \greekBo] =
\cal Q[\greekBo]$ and any $\scr B$"=element $(\greekAo, \greekBo)$
become a gauge freedom. The theorem has thus been proved. \hfill
$\blacksquare$

\vskip0.2ex
\textbf{Remark~5.} Apart from the
\hyperlink{boB}{$\brown\boB\!\!$}"=principle, the most important
aspects are: 1) \ceils{separability \tplus\ elimination} of the
$\greekA$"=symbols via cancellation of the factors
$\cal Q[\greekA\,]\hspace{0.07em}
\smash{\llap{\rotatebox{162}{$\rule{21pt}{0.5pt}$}}}$, and 2) the
fact that the existence of the function
$\cal Q[\msf c\bcdot\greekC]\;\;\forall\,\greekC \in V$ is not the same
as the existence of $\cal Q\big[\msf a\bcdot\greekA \+
\msf b\bcdot\greekB\big]$ for all $\msf a, \msf b \in\bbR$
(\hyperlink{R2}{\black Remark~\brown\textbf{2}}). Simply put, the
additive structure of the linear functional \eqref{Lf} does not bear
on `Pythagorean additivity'~\eqref{arrow} \cite{br2}. The
$\cal Q[{\cdot}{\cdot}{\cdot}]$"=functional is not linear.%
\hypertarget{I5}{}

\vskip0.1ex
\noindent%
\bfcour{\yblue Theorem $\goto$ Definitions.}
The power"=symbol notation $\|{\cdot}{\cdot}{\cdot}\|^2$ suggests now
itself, but it is not the squaring of a `something'; say, not
squaring of a norm. More to the point. The well"=definiteness and
\hyperlink{L2}{\black Lemma~\brown\textbf{2}} imply that the
$\mco C$"=functionality axiom \eqref{3} is not an assumption and can
therefore be removed from \hyperlink{D2}{\black
Definition~\brown\textbf{2}}. Also, \eqref{4*} follows from
\eqref{3}--\eqref{4}, as will be shown below. In turn, for the
$\cal Q$"=calculus to be meaningful, the
$\{\cal Q,\scr B\}$"=independence must be integral to the very
conception of calculus as a formal term. Calculus devoid invariance
and well"=definiteness \eqref{well} is no calculus at all; we should
adopt the \hyperlink{boB}{$\brown\boB\!\!$}"=invariance as a
convention. That being so, we can reformulate the Pythagorean
theorem itself.

\vskip0.2ex
\hypertarget{C1}{\textbf{Corollary~1.}} %
\emph{In the customary notation $\|\msf a\,\bcdot\,\greekA\|^2 \DEF
\cal Q[\msf a\bcdot\greekA\,]$ for the $\boB\![4]$"=invariant
$\cal Q$"=measure, the standard Pythagorean statement \eqref{pif}
ceases to be a theorem \textup{\cite{halm}} and can be replaced with a
definition of additivity \eqref{4}}:\\
\vbox to4.5ex{
\vspace{-0.7ex}
\begin{equation*}
\|[1]\msf a\bcdot\greekA \+
\msf b\bcdot\greekB\|[1]^2 = \|\msf a\bcdot\greekA\|^2 +
\|\msf b\bcdot\greekB\|^2 \qquad\forall\,\msf a,\,\msf b \in\bbR\.,
\end{equation*}}
\noindent
\emph{followed by the induced definition of perpendicular}
$\greekA\.\Over[1.8ex]{\sss\smash{\text{\!def\,}}}{\bot}\.\greekB$.

This is a rephrasing of the initial formulation \eqref{arrow};
compare with Exs.~4(a,~c) on page~123 in Ref.~\cite{halm}.

Unsurprisingly, the quantity $\msf a^2$ reframes the geometrical
concepts of a length and area merely as derivatives of the
(algebraic) square and perpendicular~$\bot$. One can also show that
the concept of an (additive) angle
$\widehat{\tgreek{a}\,\,\tgreek{g}}$ is also a derivative of the
Pythagorean square; not the other way round. As above, the right
angle $\widehat{\tgreek{a}\,\,\tgreek{g}} = 90^{\circ}$~($\Box$)
arises
\emph{prior to} the very notion of the~%
angle. 
All this gives rise to a coherent terminology on top of the
\lvs"=language: norm, distance, metric geometry, \etc. Moreover, this
makes it possible to formalize the intuitive notion of interpreting
a vector within a well"=defined framework.

Indeed, the notion of a quantity is determined by the
$\cal Q$"=scalars $\|\msf{c}\,\,{\bcdot}\,\,\greekC\+\cdots\|^2$ rather
than by symbols $\greekC$, $\|\greekC\|$, or $\msf c$. The latter
are not certain geometric `dilatators' $\hat{\msf c}$ of the
`directed segments' $\greekC$ that are `allegedly endowed' with `a
quantity'---the `length' $\|\greekC\|$; an ill"=defined language.
Returning to ``the interpretable'', the uniqueness of the measure
means that nothing but $\cal Q$ itself is to be interpreted
quantitatively. In a word, \emph{the interpretable is the $\cal Q$}:

\vbox to5.4ex{
\vspace{-0.3ex}
\begin{equation*}
\!\!\text{\small\ceils{\lvs\.\,\tplus\.\,$\scr B$}\,%
$\Under[0.9ex]{\![3]\sss\text{automorphisms}}
{\Over[1.5ex]{\![3]\sss\text{quantifying}}{\goto[10]}}$\,%
\ceils{\!$\boB\![3]$"=invariance\!}\goto[1]\ceils{interpretable%
~$\|{\cdot}{\cdot}{\cdot}\|^2$}.}
\end{equation*}}

Stated in its most essential terms, the quantification of
multi"=dimensionality is not a trivial procedure. It calls for an
intrinsic `definitio' that turns the notion into a concept.
Following \hyperlink{D1}{\black Definition~\brown$\bo1$}, we are led
to the formalization, now free of geometric motives
and imagery.%
\hypertarget{D3}{}

\smallskip
{\textbf{Corollary~2~(Definition~3)}.} %
\emph{The only well"=defined%
, meaningful and%
\ quantitative \textup{(}in the sense of
\hyperlink{D1}{\black Definition~\brown$\bo1$}\!\textup{)} %
interpretation of an abstract vector is the $\boB\![4]$"=invariant map
\eqref{xi} on orthogonal $\bo\aleph$"=bases, with the quadratic
$\mco Q$"=representation \eqref{Q} and a free
overall gauge factor $\cal Q[\greekA\,]$.}%
\hypertarget{I6}{}

\vskip0.2ex
\noindent%
\bfcour{\yblue The main theorem:~~canonicity.} The findings
presented above point out that the \lvs"=view of Pythag\-o\-rean
triangles
can, in fact, be reduced to a single act---to equip the \lvs\ with a
sole superstructure $\cal S$. Thereby, both the additivity and
\hyperlink{D2}{\black Definition\/~\brown$\bo2$} might be dispensed
with, despite their naturalness.
In keeping with the \hyperlink{boB}{$\brown\boB\!\!$}"=principle and
the topological conventions (\hyperlink{R3}{\black
Remark~\brown\textbf{3}}), we can now give \hyperlink{P}{\black
Proposition} its precise form.

\smallskip
\hypertarget{T2}{\textbf{Theorem~2.}} \itshape%
A compatible extension of the $V\![3]$"=space algebra
$\bo{\frak A}\langle V,\+,\vecO;\bbR\rangle$ by the semigroup
$\cal S$ induces, canonically, the following invariant
structures\/$:$\\[0.4ex]
\phantom{}\quad$\bullet$\;\parbox[t]{79mm}{The functional add-on
$\cal Q[{\cdot}{\cdot}{\cdot}]$ that uniquely entails the
Pythagorean squares $\|{\cdot}{\cdot}{\cdot}\|^2$ and the associated
geometry\/$:$ \textup{(}\!non\textup{)}parallelism, right triangles,
rotation, \etc.}\\[0.5ex]
\phantom{}\quad$\bullet$\;\parbox[t]{79mm}{The orthogonal bases as the
existence domain of the $\bo\aleph$"=map \eqref{xi}, followed by the
concept of a scalar product $\langle\greekA,\greekC\rangle\in\bbR$,
orthogonality
$\bot$, and the angle $\widehat{\tgreek{a}\,\,\tgreek{g}}$.}\\[0.5ex]
\phantom{}\quad$\bullet$\;\parbox[t]{79mm}{The well"=defined numeric %
\textup{(}\!quantitative\textup{)} interpretation of vectors via the
$\ds\bbR^{\!\sp}_0\!\!$"=measure $\cal Q$\/$;$
\hyperlink{D3}{\black Definition~\brown$\bo3$}.}\\[0.5ex]
This \textup{(}\!\!canonical\textup{)} theory has a unique model up to gauge
transformations $\ds \bb{GL}_2{\,\Times\.}
\bbR_{\smash{\sss\!\times}}^{\!\sp} \cong
\mrm{Aut}(V){\,\Times\,}\mrm{Aut}(\cal S)$ of the free parameters
$(\greekAo, \greekBo) \in \scr B$, $\ds\cal Q[\greekAo] \in
\bbR^{\!\sp}_0\!,$ and $\cal Q[\,\vecO\,]=0$.

\vskip0.2ex
\upshape%
We proceed as follows. Assume the $\boA$"=algebra to be compatible
with $\cal S$. Can one then lift the additive property of an
$(\greekAo,\greekBo)$"=point---an implication of
\hyperlink{L1}{\black Lemma~\brown$\bo1$}---to the level of a
structure(s) on $V$? To do this, we seek these extensions within the
class of canonical/\!"!natural ones, since their unambiguous link with
the original data \ceils{$V\!\!$ \tplus\ $\cal S$} guarantees
functoriality of the whole derivation.

\emph{Proof} of \hyperlink{T2}{\black Theorem~\brown$\bo2$}. Let us
augment the signature of the $\boA$"=algebra with a functional symbol
$\cal Q$ of the $\cal S$"=superstructure. The resulting algebra
$\ds\boAo\!\big\langle\!\!(V,\+,\vecO;\bbR)$, $\ds(\cal
Q{\.\in\.}\bbR^{\!\sp}_0\!;+)\!\!\big\rangle$ also incorporates an axiom
with variables $\{\greekAo,\greekBo\}$ bound by an existential
quantifier $\exists (\greekAo,\greekBo){:}\;\cal
Q[\greekAo\+\greekBo]= \cal Q[\greekAo] + \cal Q[\greekBo]$ and will
be regarded as canonical.\@\ To extend this $\scr A\!$"=property, the
canonicity requirement suggests a natural passage from a certain
point $(\greekAo,\greekBo)\in V^2$ to the canonical domain for the
algebra $\boAo\!\langle{\cdot}{\cdot}{\cdot},\VV\rangle$---the whole
space $(\VV)$. This is to supplement this $\exists$"=axiom with the
definition of a function:

\vbox to4.6ex{
\begin{equation*}
V^2\mapsto\bbR{:}\;\;\Delta(\greekA,\greekB)\DEF
\cal Q[\greekA\+\greekB]-\cal Q[\greekA\,]-\cal Q[\greekB]\;\;\;\forall\,
\greekA,\greekB\in V.
\end{equation*}}
\vskip0.5ex
\noindent
The $\scr A\!$"=pair corresponds to the zero level
$\Delta(\greekAo,\greekBo) = 0_{\raisebox{0.1ex}{$\sss\bbR$}}$.

Clearly, $\Delta$ is a canonical extension of
$\boAo\!\langle{\cdot}{\cdot}{\cdot},\VV\rangle$, which has the
structural automorphisms $V^2\.{\mapsto}\. V^2$ in their own right.
They are induced by the canonical automorphisms of
$\boA\langle{\cdot}{\cdot}{\cdot}\rangle$, yielding
$(\greekA,\greekB)\to\hat{\text{\small$\msf A$}}(\greekA,\greekB)
\mathbin{\Over{\text{def}}{=}} (\hat{\msf a}\,\greekA$,
$\skew{-3}\hat{\msf b}\,\greekB)$, and by the permutation on the
$(\VV)$"=set: $\hat{\text{\small$\frak T$}}(\greekA,\greekB) =
(\greekB,\greekA)$.

It is also clear that in category $\mathbf{Vect}
_{\raisebox{0.1ex}{$\sss\bbR,\!\cal Q$}}$, with objects $(V,\cal Q)$
and linear morphisms respecting the $\cal Q$"=fibres \eqref{well},
$\Delta$ is a canonical morphism; $\cal Q$ is not. Since such
morphisms commute with ones from terminal objects, the
signature"=constant status is preserved. In particular, the
`zero'"=constant $(\vecO,\vecO)$ of the algebra
$\boA_{\sss\Delta}\![3]\langle \VV,
\Delta,{\cdot}{\cdot}{\cdot}\rangle$ is mapped to the $\bbR$"=zero:
$\Delta(\vecO,\vecO)=0_{\raisebox{0.1ex}{$\sss\bbR$}}$. On the other
hand, $\Delta(\vecO,\vecO)=-\cal Q[\,\vecO\,]$, whence the canonical
correspondence between the `zeroes':
$\cal Q[\,\vecO\,]=0_{\raisebox{0.1ex}{$\sss\bbR$}}$. Thus, in
addition to $(\greekAo,\greekBo)$, we obtain another point,
$(\vecO,\vecO)$, of the level~$\Delta=0$.

Consider now automorphisms of $\boA$ and $\boA_{\sss\Delta}$.
Compatibility \eqref{well} and \hyperlink{L2}{\black
Lemma~\brown\textbf{2}} lead to axiom~\eqref{3}. Obviously,
$\Delta\circ\hat{\text{\small$\frak T$}}\equiv\Delta$. As for
$\hat{\text{\small$\msf A$}}$, canonicity implies
\ceils{$\Delta(\greekAo,\greekBo) = 0 = \Delta(\vecO,\vecO)$}
\hence\ \ceils{$(\Delta\circ\hat{\text{\small$\msf
A$}})(\greekAo,\greekBo) = (\Delta\circ\hat{\text{\small$\msf
A$}})(\vecO,\vecO)$} \hence\
$\big|(\Delta\circ\hat{\text{\small$\msf A$}})(\vecO,\vecO) =
0\big|$ \hence\ \ceils{$\Delta(\msf a\bcdot\greekAo, \msf
b\bcdot\greekBo)=0\;\.\forall\,\msf a,\msf b$}. This is axiom
\eqref{4}. Conditions \eqref{3}--\eqref{4} and \hyperlink{I}{\black
Theorem~\brown\textbf{1}\cour{{\![3](\!I\!)}}} entail that
$\greekBo\ne\msf k\bcdot\greekAo$; call this basis. The function
$\cal Q[\greekC]=\bo\aleph(\msf a, \msf b; \greekAo, \greekBo)$ has
thus been defined for all $\greekC\in V\!\!$, in accordance with the
$\boAo\!$"=signature for $\cal Q$: $\ds V\mapsto\bbR^{\!\sp}_0$.
Moreover, any of the bases $(\greekA',\greekB')=(\msf
a\bcdot\greekAo, \msf b\bcdot\greekBo)\in\scr B$ may serve as a seed
$\scr A\!$"=pair;
we drop subscripts further.

Although the scalability
$\cal Q[\greekA']=\msf a^{2\textsf{p}}{\,\times}\cal Q[\greekA\,]$
introduces the $\smallc({\times}\smallc)$"=multiplication
(\hyperlink{I}{\black Theorem~\brown\textbf{1}\cour{{\![3](\!I\!)}}}\![4]),
the structure $\ds\langle\, \bbR^{\!\sp}\!, \times\rangle\FED
\bbR_{\smash{\sss\!\times}}^{\!\sp}$ here is canonically isomorphic to
the group $\mrm{Aut}(\cal S)$. Then, the original
$(\greekA,\greekB)$"=basis may be renormalized to an $\scr A\!$"=pair
$\{\greekA', \greekB'\}$ with any, \eg, equal $\cal Q$"=magnitudes.
Indeed, $\{\cal Q[\greekA']$ $=$ $\cal Q[\greekB']\}$ \hhence\
$\{\msf a^{2\textsf{p}}{\,\cdot\,}\cal Q[\greekA\,] =
\msf b^{\textsf{2p}}{\,\cdot\,}\cal Q[\greekB]\}$. From this, numbers
$(\msf a,\msf b)$ can always be found for any
$(\cal Q[\greekA\,],\cal Q[\greekB])$ by the monotonicity of the
function $\msf a^{2\textsf{p}}$ on the interval
$\msf a^2\in[0,\infty)$. Taking that $(\greekA',\greekB')$ as an
$\scr A\!$"=pair, we derive: $\cal Q\big[\msf a\bcdot\greekA' \+
\msf b\bcdot\greekB'\big] = \big| \eqref{4}\big| =
\msf a^{2\textsf{p}}{\,\cdot\,}\cal Q[\greekA'] +
\msf b^{2\textsf{p}}{\,\cdot\,}\cal Q[\greekB'] = \big|\cal Q[\greekA']
\rightleftarrows \cal Q[\greekB'] \big| =
\msf a^{2\textsf{p}}{\,\cdot\,}\cal Q[\greekB'] +
\msf b^{2\textsf{p}}{\,\cdot\,}\cal Q[\greekA'] =
\cal Q[\msf a\bcdot\greekB'] + \cal Q[\msf b\bcdot\greekA'] = \big|
\eqref{4}\big| = \cal Q\big[\msf a\bcdot\greekB' \+
\msf b\bcdot\greekA'\big] = \cal Q\big[\msf a\bcdot\greekA' \+
\msf b\bcdot\greekB'\big]$. Up to the prime"=symbol, this is axiom
\eqref{4*}; in fact, a consequence of \eqref{3}--\eqref{4}.
\hyperlink{D2}{\black Definition~\brown\textbf{2}} is now fully
deduced.

The existence of a model for these canonical constructions is
ensured by the \hyperlink{boB}{$\brown\boB\!\!$}"=invariance and
\hyperlink{II}{\black
Theorem~\brown\textbf{1}\cour{{\![3](\!I\!I\!}}}--%
\hyperlink{III}{\black \brown\cour{{\!I\!I\!I\!)}}}. Canonicity
introduces no arbitrary elements, and so the models have no other
degrees of freedom beyond the $\exists$"=axiom about
$(\greekAo,\greekBo)$ and $\cal Q[\greekAo]$. The corresponding
gauge transformations are thus the $\bb{GL}_2$"=change of bases and
the $\ds\bbR_{\smash{\sss\!\times}}^{\!\sp}$"=multiplication
$\cal Q[\greekC] \goto \cal Q[\greekAo] {\,\times\,} \cal Q[\greekC]$.

In summary, the $\Delta$"=structure describes an axiom"=free
interaction of $\cal S$ with algebra $\boA\langle
V,{\cdot}{\cdot}{\cdot}\rangle$ and, depending on context, may be
thought of as a deviation $\Delta$ from additivity
$\cal Q[\greekA\+\greekC] = (\cal Q[\greekA\,] + \cal Q[\greekC]) +
\Delta(\greekA,\greekC)$ or referred to as the scalar product
$\langle\greekA,\greekC\rangle \DEF \cal Q[\greekA\+\greekC] -
\cal Q[\greekA\,]-\cal Q[\greekC]$. It is positive"=definite and
results in the formal orthogonality $(\greekA\.\bot\.\greekC)\hhence
\langle\greekA,\greekC\rangle=0$, the concept of an angle
$\widehat{\tgreek{a}\,\,\tgreek{g}}$, rotation, \etc. This also
accounts for the (canonical) origin of the $\cal Q$"=polarization and
the parallelogram law, and thereby answers the question posed in the
\hyperlink{I1}{\black Introduction}: these relations/\!"!terms are not
necessary for the theory. \hfill $\blacksquare$

\smallskip
\hypertarget{C3}{\textbf{Corollary~3.}} %
\itshape Mutatis mutandis, both the theory and \hyperlink{D3}{\black
Definition~\,\brown$\bo3$} extend to the \textup{(}\![3]`quantum\!'\textup{)}
field\/~$\bbC$ and to the finite dimensions ${\dim}\,\, V\![3]$, thereby
yielding the Born rule.

\upshape An additional technical point here is the extra
\emph{canonical} automorphism of the algebra $\boA\langle
V,\+,\vecO;\bbC\rangle$---the complex conjugation
$\ds\bbC\mapsto\bbC^*$. This yields $\msf a^2\goto|\frak a|^2$
\cite{br2}.

\hypertarget{I7}{}

\end{document}